


\magnification 1200
\hsize 13.2cm
\vsize 20cm
\parskip 3pt plus 1pt
\parindent 5mm

\def\\{\hfil\break}


\font\seventeenbf=cmbx10 at 17.28pt

\font\twelvebf=cmbx10 at 12pt
\font\eightbf=cmbx8
\font\sixbf=cmbx6

\font\eighti=cmmi8
\font\sixi=cmmi6

\font\eightrm=cmr8
\font\sixrm=cmr6

\font\eightsy=cmsy8
\font\sixsy=cmsy6

\font\eightit=cmti8
\font\eighttt=cmtt8
\font\eightsl=cmsl8

\font\seventeenbsy=cmbsy10 at 17.28pt

\font\twelvebsy=cmbsy10 at 12pt
\font\tenbsy=cmbsy10
\font\eightbsy=cmbsy8
\font\sevenbsy=cmbsy7
\font\sixbsy=cmbsy6
\font\fivebsy=cmbsy5

\font\tenmsa=msam10

\font\sevenmsa=msam7
\font\fivemsa=msam5
\newfam\msafam
  \textfont\msafam=\tenmsa
  \scriptfont\msafam=\sevenmsa
  \scriptscriptfont\msafam=\fivemsa

\font\tenmsb=msbm10
\font\eightmsb=msbm8
\font\sevenmsb=msbm7
\font\fivemsb=msbm5
\newfam\msbfam
  \textfont\msbfam=\tenmsb
  \scriptfont\msbfam=\sevenmsb
  \scriptscriptfont\msbfam=\fivemsb
\def\Bbb{\fam\msbfam\tenmsb}

\font\tenCal=eusm10
\font\sevenCal=eusm7
\font\fiveCal=eusm5
\newfam\Calfam
  \textfont\Calfam=\tenCal
  \scriptfont\Calfam=\sevenCal
  \scriptscriptfont\Calfam=\fiveCal
\def\Cal{\fam\Calfam\tenCal}

\font\teneuf=eusm10
\font\teneuf=eufm10
\font\seveneuf=eufm7
\font\fiveeuf=eufm5
\newfam\euffam
  \textfont\euffam=\teneuf
  \scriptfont\euffam=\seveneuf
  \scriptscriptfont\euffam=\fiveeuf

\font\seventeenbfit=cmmib10 at 17.28pt

\font\twelvebfit=cmmib10 at 12pt
\font\tenbfit=cmmib10
\font\eightbfit=cmmib8
\font\sevenbfit=cmmib7
\font\sixbfit=cmmib6
\font\fivebfit=cmmib5
\newfam\bfitfam
  \textfont\bfitfam=\tenbfit
  \scriptfont\bfitfam=\sevenbfit
  \scriptscriptfont\bfitfam=\fivebfit


\catcode`\@=11
\def\eightpoint{%
  \textfont0=\eightrm \scriptfont0=\sixrm \scriptscriptfont0=\fiverm
  \def\rm{\fam\z@\eightrm}%
  \textfont1=\eighti \scriptfont1=\sixi \scriptscriptfont1=\fivei
  \def\oldstyle{\fam\@ne\eighti}%
  \textfont2=\eightsy \scriptfont2=\sixsy \scriptscriptfont2=\fivesy
  \textfont\itfam=\eightit
  \def\it{\fam\itfam\eightit}%
  \textfont\slfam=\eightsl
  \def\sl{\fam\slfam\eightsl}%
  \textfont\bffam=\eightbf \scriptfont\bffam=\sixbf
  \scriptscriptfont\bffam=\fivebf
  \def\bf{\fam\bffam\eightbf}%
  \textfont\ttfam=\eighttt
  \def\tt{\fam\ttfam\eighttt}%
  \textfont\msbfam=\eightmsb
  \def\Bbb{\fam\msbfam\eightmsb}%
  \abovedisplayskip=9pt plus 2pt minus 6pt
  \abovedisplayshortskip=0pt plus 2pt
  \belowdisplayskip=9pt plus 2pt minus 6pt
  \belowdisplayshortskip=5pt plus 2pt minus 3pt
  \smallskipamount=2pt plus 1pt minus 1pt
  \medskipamount=4pt plus 2pt minus 1pt
  \bigskipamount=9pt plus 3pt minus 3pt
  \normalbaselineskip=9pt
  \setbox\strutbox=\hbox{\vrule height7pt depth2pt width0pt}%
  \let\bigf@ntpc=\eightrm \let\smallf@ntpc=\sixrm
  \normalbaselines\rm}
\catcode`\@=12

\def\eightpointbf{%
 \textfont0=\eightbf   \scriptfont0=\sixbf   \scriptscriptfont0=\fivebf
 \textfont1=\eightbfit \scriptfont1=\sixbfit \scriptscriptfont1=\fivebfit
 \textfont2=\eightbsy  \scriptfont2=\sixbsy  \scriptscriptfont2=\fivebsy
 \eightbf
 \baselineskip=10pt}

\def\tenpointbf{%
 \textfont0=\tenbf   \scriptfont0=\sevenbf   \scriptscriptfont0=\fivebf
 \textfont1=\tenbfit \scriptfont1=\sevenbfit \scriptscriptfont1=\fivebfit
 \textfont2=\tenbsy  \scriptfont2=\sevenbsy  \scriptscriptfont2=\fivebsy
 \tenbf}

\def\twelvepointbf{%
 \textfont0=\twelvebf   \scriptfont0=\eightbf   \scriptscriptfont0=\sixbf
 \textfont1=\twelvebfit \scriptfont1=\eightbfit \scriptscriptfont1=\sixbfit
 \textfont2=\twelvebsy  \scriptfont2=\eightbsy  \scriptscriptfont2=\sixbsy
 \twelvebf
 \baselineskip=14.4pt}

\def\seventeenpointbf{%
 \textfont0=\seventeenbf  \scriptfont0=\twelvebf  \scriptscriptfont0=\eightbf
 \textfont1=\seventeenbfit\scriptfont1=\twelvebfit\scriptscriptfont1=\eightbfit
 \textfont2=\seventeenbsy \scriptfont2=\twelvebsy \scriptscriptfont2=\eightbsy
 \seventeenbf
 \baselineskip=20.736pt}


\newdimen\srdim \srdim=\hsize
\newdimen\irdim \irdim=\hsize
\def\NOSECTREF#1{\noindent\hbox to \srdim{\null\dotfill ???(#1)}}
\def\SECTREF#1{\noindent\hbox to \srdim{\csname REF\romannumeral#1\endcsname}}
\def\INDREF#1{\noindent\hbox to \irdim{\csname IND\romannumeral#1\endcsname}}
\newlinechar=`\^^J
\def\openauxfile{
  \immediate\openin1\jobname.aux
  \ifeof1
  \message{^^JCAUTION\string: you MUST run TeX a second time^^J}
  \let\sectref=\NOSECTREF \let\indref=\NOSECTREF
  \else
  \input \jobname.aux
  \message{^^JCAUTION\string: if the file has just been modified you may
    have to run TeX twice^^J}
  \let\sectref=\SECTREF \let\indref=\INDREF
  \fi
  \message{to get correct page numbers displayed in Contents or Index
    Tables^^J}
  \immediate\openout1=\jobname.aux
  \let\END=\end \def\end{\immediate\closeout1\END}}

\newbox\titlebox   \setbox\titlebox\hbox{\hfil}
\newbox\sectionbox \setbox\sectionbox\hbox{\hfil}
\def\folio{\ifnum\pageno=1 \hfil \else \ifodd\pageno
           \hfil {\eightpoint\copy\sectionbox\kern8mm\number\pageno}\else
           {\eightpoint\number\pageno\kern8mm\copy\titlebox}\hfil \fi\fi}
\footline={\hfil}
\headline={\folio}

\def\titlerunning#1{\setbox\titlebox\hbox{\eightpoint #1}}
\def\title#1{\noindent\hfil$\smash{\hbox{\seventeenpointbf #1}}$\hfil
             \titlerunning{#1}\medskip}

\newcount\numbersection \numbersection=-1
\def\sectionrunning#1{\setbox\sectionbox\hbox{\eightpoint #1}
  \immediate\write1{\string\def \string\REF
      \romannumeral\numbersection \string{%
      \noexpand#1 \string\dotfill \space \number\pageno \string}}}
\def\section#1{%
  \par\vskip0.666cm\penalty -100
  \vbox{\baselineskip=14.4pt\noindent{{\twelvepointbf #1}}}
  \vskip2pt
  \penalty 500
  \advance\numbersection by 1
  \sectionrunning{#1}}

\def\subsection#1{%
  \par\vskip0.5cm\penalty -100
  \vbox{\noindent{{\tenpointbf #1}}}
  \vskip1pt
  \penalty 500}

\newcount\numberindex \numberindex=0
\def\index#1#2{%
  \advance\numberindex by 1
  \immediate\write1{\string\def \string\IND #1%
     \romannumeral\numberindex \string{%
     \noexpand#2 \string\dotfill \space \string\S \number\numbersection,
     p.\string\ \space\number\pageno \string}}}

\newdimen\itemindent \itemindent=\parindent

\def\item#1{\par\noindent\hangindent\itemindent%
            \rlap{#1}\kern\itemindent\ignorespaces}
\def\itemitem#1{\par\noindent\hangindent2\itemindent%
            \kern\itemindent\rlap{#1}\kern\itemindent\ignorespaces}
\def\itemitemitem#1{\par\noindent\hangindent3\itemindent%
            \kern2\itemindent\rlap{#1}\kern\itemindent\ignorespaces}

\long\def\claim#1|#2\endclaim{\par\vskip 5pt\noindent
{\tenpointbf #1.}\ {\it #2}\par\vskip 5pt}

\def\proof{\noindent{\it Proof}}

\def\today{\ifcase\month\or
January\or February\or March\or April\or May\or June\or July\or August\or
September\or October\or November\or December\fi \space\number\day,
\number\year}

\catcode`\@=11
\newcount\@tempcnta \newcount\@tempcntb
\def\timeofday{{%
\@tempcnta=\time \divide\@tempcnta by 60 \@tempcntb=\@tempcnta
\multiply\@tempcntb by -60 \advance\@tempcntb by \time
\ifnum\@tempcntb > 9 \number\@tempcnta:\number\@tempcntb
  \else\number\@tempcnta:0\number\@tempcntb\fi}}
\catcode`\@=12

\def\bibitem#1&#2&#3&#4&%
{\hangindent=0.8cm\hangafter=1
\noindent\rlap{\hbox{\eightpointbf #1}}\kern0.8cm{\rm #2}{\it #3}{\rm #4.}}


\def\bC{{\Bbb C}}

\def\bQ{{\Bbb Q}}


\def\cO{{\Cal O}}

\def\cV{{\Cal V}}


\def\square{{\hfill \hbox{
\vrule height 1.453ex  width 0.093ex  depth 0ex
\vrule height 1.5ex  width 1.3ex  depth -1.407ex\kern-0.1ex
\vrule height 1.453ex  width 0.093ex  depth 0ex\kern-1.35ex
\vrule height 0.093ex  width 1.3ex  depth 0ex}}}
\def\qed{\kern10pt$\square$}
\def\hexnbr#1{\ifnum#1<10 \number#1\else
 \ifnum#1=10 A\else\ifnum#1=11 B\else\ifnum#1=12 C\else
 \ifnum#1=13 D\else\ifnum#1=14 E\else\ifnum#1=15 F\fi\fi\fi\fi\fi\fi\fi}
\def\msatype{\hexnbr\msafam}
\def\msbtype{\hexnbr\msbfam}
\mathchardef\restriction="3\msatype16   
\mathchardef\boxsquare="3\msatype03
\mathchardef\preccurlyeq="3\msatype34
\mathchardef\compact="3\msatype62
\mathchardef\smallsetminus="2\msbtype72   
\mathchardef\subsetneq="3\msbtype28
\mathchardef\supsetneq="3\msbtype29
\mathchardef\leqslant="3\msatype36   
\mathchardef\geqslant="3\msatype3E   
\mathchardef\stimes="2\msatype02
\mathchardef\ltimes="2\msbtype6E
\mathchardef\rtimes="2\msbtype6F

\def\dbar{\overline\partial}
\def\ddbar{\partial\overline\partial}

\let\ol=\overline

\let\wt=\widetilde

\let\text=\hbox
\def\buildo#1^#2{\mathop{#1}\limits^{#2}}
\def\buildu#1_#2{\mathop{#1}\limits_{#2}}
\def\ort{\mathop{\hbox{\kern1pt\vrule width4.0pt height0.4pt depth0pt
    \vrule width0.4pt height6.0pt depth0pt\kern3.5pt}}}

\def\vlra{\mathrel{\smash-}\joinrel\mathrel{\smash-}\joinrel%
    \kern-2pt\longrightarrow}
\def\srelbar{\vrule width0.6ex height0.65ex depth-0.55ex}
\def\merto{\mathrel{\srelbar\kern1.3pt\srelbar\kern1.3pt\srelbar
    \kern1.3pt\srelbar\kern-1ex\raise0.28ex\hbox{${\scriptscriptstyle>}$}}}



\def\Supp{\mathop{\rm Supp}\nolimits}


\long\def\InsertFig#1 #2 #3 #4\EndFig{\par
\hbox{\hskip #1mm$\vbox to#2mm{\vfil\special{"
(/home/demailly/psinputs/grlib.ps) run
#3}}#4$}}
\long\def\LabelTeX#1 #2 #3\ELTX{\rlap{\kern#1mm\raise#2mm\hbox{#3}}}


\itemindent = 7mm

\title{Bergman kernels and}
\title{subadjunction}
\titlerunning{Subadjunction}

\vskip10pt

\centerline {\tenrm Bo BERNDTSSON and Mihai P\u AUN}
á



\vskip20pt

\noindent{\bf Abstract. \it {In this article our main result is a more complete version of the statements obtained in {\rm [6]}. One of the important technical point of our proof is an 
$\displaystyle L^{2\over m}$ extension theorem of Ohsawa-Takegoshi type, which is derived from the original result by a simple fixed point  method.
Moreover, we show that these techniques combined with an appropriate form of the``invariance of plurigenera" can be used in order to obtain a new proof of the celebrated
Y. Kawamata subadjunction theorem.
}}

%
\section{\S0 Introduction}

\medskip
\noindent Our starting point in this article is the following
generalization of our work [6].

\claim 0.1 Theorem|Let $p: X\to Y$ be a surjective projective map between two non-singular manifolds, and let 
$L\to X$ be a line bundle endowed with a metric $h_L$ such that:

{\itemindent 4mm
\smallskip
\item {\rm 1)} The curvature current of the bundle $(L, h_L)$ is semipositive, i.e. 
$\displaystyle \sqrt {-1}\Theta_{h_L}(L)\geq 0$;

\item {\rm 2)} There exists a generic point $z\in Y$ and a section 
$\displaystyle u\in H^0\big(X_z, mK_{X_z}+ L\big)$
such that 
$$\displaystyle \int_{X_z}|u|^{2\over m}e^{-{{\varphi_L}\over {m}}}< \infty.$$

}

\noindent Then the line bundle $mK_{X/Y}+ L$ admits a metric with positive curvature current.
Moreover, this metric is equal to the fiberwise $m$--Bergman kernel metric 
on the generic fibers of $p$. \hfill\qed
\endclaim
\noindent In the statement above the meaning of the word {\sl generic} is as follows: $w\in Y$ is generic if it is not a critical value of $p$ and if
the sections of the bundle $mK_{X/Y}+ L_{|X_w}$ extend locally near $w$. 
\smallskip

\noindent This kind of positivity results (and their relevance for important problems in algebraic geometry) have been investigated since quite long time; we will only mention here a few contributors [8], [13], [14], [15], [16], [17], [19], [20], [21], [22], [23], [26], [27], [28], [31],
[34], [35] (and we apologize to the ones we omit...). For example, one can see that if the metric of $L$ is given by a section of some of its multiples, then the  
{\sl qualitative part} of the theorem above can be obtained as a consequence of 
E. Viehweg's results on weak positivity of direct images (cf. [35] and 
the references therein, as well as the recent results of A. H\"oring [16]).  

\noindent The difference between theorem 0.1 and the corresponding result we obtain in [6] is that 
here the map $p$ can be singular (in our previous work we have assumed that $p$ is a smooth fibration). 
Actually, the additional technical result allowing us to treat the general case
is an $L^{2/m}$ extension theorem which we describe next. 

The set-up is the following: let $\Omega\subset \bC^n$ be a ball of radius $r$ and let $h:\Omega\to \bC$ be a holomorphic function, such that 
$\sup_\Omega|h|\leq 1$; moreover, we assume that the gradient $\partial h$ of $h$ 
is nowhere zero on the set $V:=  (h=0)$. We denote by $\varphi$ a plurisubharmonic function, such that its restriction to $V$ is well-defined (i.e., $\varphi_{|V}\neq -\infty$). Then the Ohsawa-Takegoshi extension theorem states that for any $f$ holomorphic on $V$, there exists a function $F$, holomorphic in all of $\Omega$, such that $F=f$ on $V$, and moreover

$$\int_{\Omega} |F|^{2}\exp (-\varphi)d\lambda \leq C_0\int_{V} |f|^{2}\exp (-\varphi){{d\lambda_V}\over {|\partial h|^2}}. $$
Here, $C_0$ is an absolute constant. Our generalization is the following.

\claim 0.2 Proposition|For any holomorphic 
function
$f: V\to \bC$  with the property that 
$$\int_{V} |f|^{2/m}\exp (-\varphi){{d\lambda_V}\over {|\partial h|^2}}\leq 1, $$
there exists a function $F\in {\cal O}(\Omega)$ such that:
\item{(i)} $F_{|V}= f$ i.e. the function $F$ is an extension of $f$;
\smallskip
\item{(ii)} The next $L^{2/m}$ bound holds 
$$\int_{\Omega} |F|^{2/m}\exp (-\varphi)d\lambda \leq C_0\int_{V} |f|^{2/m}\exp (-\varphi){{d\lambda_V}\over {|\partial h|^2}}, $$
where $C_0$ is the same constant as in the Ohsawa-Takegoshi theorem.
\endclaim

Once this result is established, the proof of 0.1 runs as follows: we first restrict ourselves to the Zariski open set $Y_0\subset Y$ which does not contain any critical value of $p$, and such that 
$\forall z\in Y_0$, all the sections of the bundle $\displaystyle mK_{X_z}+ L$ extend locally near $z$ (the existence of such a set $Y_0$ is a consequence of standard semi-continuity 
arguments). Over $Y_0$, we can apply the results in [6] and therefore the $m$--Bergman kernel metric has a psh variation. We then use the  $L^{2/m}$ extension result to estimate our metric from above by a uniform constant. Standard results of pluripotential theory then gives that the metric extends to a semipositively curved  metric on all of $X$. At the end of this section, we derive some corollaries and further results; we equally provide 
some additional explanations concerning our previous work [6].
\hfill\qed

\medskip 
\noindent  We show next that the techniques and results we obtain here can be used in order to provide a direct and transparent argument for the next statement 0.3,
which appears to be the main part of the proof of the adjunction result originally due to Y. Kawamata (cf. [19]). We also refer to the recent presentation of J. Koll\'ar (see [23]) and the references therein, especially the recent developments due to F. Ambro and O. Fujino. We will adopt the terminology in [19], [23], so that the relation between our arguments and the construction of the proof in these articles becomes as clear as possible.

Let $X$ be a {\sl normal projective variety}, and let $\Delta$ be an effective Weil $\bQ$-divisor
on $X$, such that $K_X+ \Delta$ is $\bQ$-Cartier. We assume {\sl for the moment} that the pair 
$(X, \Delta)$ is {\sl log-canonical}. This requirement means that 
there exists a log-resolution $\mu: X^\prime\to X$ of the pair $(X, \Delta)$, such that
we have
$$\mu^\star(K_X+ \Delta)= K_{X^\prime}+ \sum_{j\in J}a^jY_j\leqno(1)$$
where the coefficients $(a^j)$ above are rational numbers, such $\forall j\in J$ we have $a^j\leq 1$. We take this opportunity to recall also that 
if the inequality $a^j< 1$ is satisfied for all $j\in J$, then the pair $(X, \Delta)$ is called {\sl Kawamata log terminal, or klt for short}. The above properties are independent of the log-resolution
$\mu$ of the pair $(X, \Delta)$ (cf [23]). 

A subvariety $W\subset X$ is called a {\sl log canonical center} if there exists 
a log resolution $\mu$ of the pair $(X, \Delta)$ and a hypersurface say $Y_1$ whose corresponding coefficient $a^1$ is equal to 1, such that $W= \mu(Y_1)$. Since each component of the intersection of two log canonical centers is a log canonical center (cf. [18], [23]), the notion of {\sl minimal log canonical} center makes sense. 

\noindent In this framework, the following important result was established by Y. Kawamata in [19].

\claim Theorem ([19])|Let $(X, \Delta)$ be a log-canonical pair; we assume moreover the existence of an effective divisor $\Delta_0\leq \Delta$, such that the pair $(X, \Delta_0)$
is klt. Let $H$ be an ample divisor  on $X$, and let $\varepsilon > 0$ be a positive rational.  
We consider a minimal log-canonical center $W$ of the pair $(X, \Delta)$. 
\\ Then $W$ is normal, and there exists an effective $\bQ$-divisor $\Delta_W$ on 
$W$ such that:

{\itemindent 9pt

\item {\rm i)} We have the adjunction relation
$K_X+ \Delta+ \varepsilon H_{|W}\equiv K_W+ \Delta_W;$

\item{\rm ii)} The pair $(W, \Delta_W)$ is klt.
\hfill\qed }
\endclaim
We can distinguish two parts in the original proof of this result: one starts with a few important reductions, showing that under the hypothesis of the preceding theorem,
$W$ becomes an {\sl exceptional minimal center} for some pair $(X, \wt\Delta)$
in the terminology used in [23]. By definition, this means that in the decomposition 
above corresponding to $K_X+ \wt\Delta$ we will have $Y_j\cap Y_1=\emptyset$, as soon as $a^j= 1$; the technique used is the so-called {\sl tie-break} method (cf. [17] [18], [23]). Also, the minimal
center $W$ is proved to be {\sl normal}, as a consequence of the Kawamata-Viehweg vanishing theorem. The hypothesis $(X, \Delta_0)$
{\sl is klt} is crucial for the proof of these results; we refer to the article by F. Ambro [1] for a generalization in the context of {\sl log-varieties} (see also the work of O. Fujino [13]), and the references therein.

After these reductions, the heart of the matter is the construction of the 
divisor $\Delta_W$. This is performed by using {\sl ramified cover tricks} (in order to deal with $\bQ$-divisors instead of line bundles), and very subtle properties derived from the existence of mixed Hodge structures, together with the work of P. Griffiths, W. Schmid, and many other contributors. 

\smallskip
\noindent In the second part of the present article, our aim is to present a {\sl metric approach} 
for the construction of the divisor $\Delta_W$ above. Therefore, we will assume from the very beginning that $W$ is an {\sl exceptional log canonical center} of $(X, \Delta)$: starting from this, we will use the theorem 0.1 above in order to obtain $\Delta_W$.
We remark that no 
ramified covers or mixed Hodge structures are used in our arguments. Our proof 
is explicit and natural in the context of Bergman kernels, 
but on the other hand, the original argument provides a finer analysis of the singularities of $\Delta_W$ on a birational model of $W$.

Actually, we will work in the following {\sl local context}. We assume that the pair
$(X, \Delta)$ is log canonical 
at some point $x_0\in X$, that is to say that 
there exists a log-resolution $\mu: X^\prime\to X$ of the pair $(X, \Delta)$, such that
we have
$$\mu^\star(K_X+ \Delta)= K_{X^\prime}+ \sum_{j\in J}a^jY_j\leqno (2)$$
where the coefficients $(a^j)$ above are rational numbers, such that $a^j\leq 1$
if $x_0\in \mu(Y_j)$. As before, the above property is independent of the log-resolution
$\mu$ of the pair $(X, \Delta)$.

By the preceding discussion, there is no loss of generality if we only focus on the properties of the {\sl exceptional centers} $W$, so we assume that
there exists a log-resolution $\mu: X^\prime\to X$ of the pair $(X, \Delta)$ together with a decomposition of the inverse image of the $\bQ$-divisor $K_X+ \Delta$ as follows
$$\mu^\star(K_X+ \Delta)=  K_{X^\prime}+ S+ \Delta^\prime+ R- \Xi,\leqno (\dagger)$$
such that:
{
\itemindent 3.5mm
\smallskip
\item {$\bullet$} $S$ is an irreducible hypersurface, such that $W= \mu(S)$;
\smallskip
\item {$\bullet$} $\Delta^\prime:= \sum_j a^jY_j$, where $x_0\in \mu(Y_j)$ and $a^j\in ]0, 1[$;
\smallskip
\item {$\bullet$} The divisor $R$ is effective, and a hypersurface $Y_j$ belongs to 
its support if either $x_0\not\in \mu(Y_j)$, or $Y_j\cap S= \emptyset$ (so in particular
the restriction $R_{|S}$ is $\mu_{|S}$-vertical);

\smallskip
\item {$\bullet$} The divisor $\Xi$ is effective and $\mu$-contractible; in addition, we assume that the divisors appearing in the right hand side of the formula $(\dagger)$ 
have strictly normal crossings.

}

\noindent In general, the center $W$ is singular, and we will assume that the restriction map
$\displaystyle \mu_{|S}: S\to W$
factors through the desingularization $g: W^\prime\to W$, so that we have
$$\mu_{|S}= g\circ p$$
where $p: S\to W^\prime$ is a surjective projective map. 

Before stating our next result, we introduce a last piece of notation: let $A$ be an ample bundle on $X$, and let $F_1,..., F_k$ be a set of smooth hypersurfaces of $W^\prime$ with strictly normal crossings, such that there exists positives rational numbers 
$(\delta_j)$ for which the $\bQ$-bundle
$$g^\star(A)-\sum \delta^jF_j+ \varepsilon K_{W^\prime}\leqno (3)$$
is semi-positive (in metric sense) for any $\varepsilon$ small enough, 
and such that $g(F_j)\subset W_{\rm sing}$ for each $j$.
Indeed a set $(F_j)$ with the properties specified above does exists, as the variety $X$ sits inside some projective space, and we can consider an embedded resolution of singularities of the set $W$. It is routine to verify that the 
intersection of the exceptional divisors with $W^\prime$ will do. We remark that if the canonical bundle $K_{W^\prime}$ happens to be the inverse image of some bundle on 
$W$, then we can take $\delta^j= 0$ for all $j$.

The family $(F_j)$ induces a decomposition of the divisor $\Xi$ as follows
$$\Xi= \Xi_1+ \Xi_2\leqno (\dagger\dagger)$$
where by definition $\Xi_{1}$ is the part of the divisor $\Xi$
whose support restricted to $S$ is  
mapped 
by $p$ into $\cup_j F_j$.

\medskip
\noindent Our version of the result of Y. Kawamata is the following.

\claim 0.3 Theorem|Under the hypotheses above,
there exists a closed positive current $\Theta_{W^\prime}$  
and an effective divisor $\Lambda$ on $W^\prime$ such that:
{\itemindent 4.5mm
\smallskip
\item {\rm (a)} We have $K_{W^\prime}+ \Theta_{W^\prime}\equiv \Lambda+ g^\star ({K_X+ \Delta})_{|W^\prime}$;
\smallskip
\item {\rm (b)} The (local) function
 $$\exp\big((1+ \varepsilon_0)(\varphi_{\Lambda}-\varphi_{\Theta_{W^\prime}})\big)$$ 
 is in $L^1_{\rm loc}$ each point of $W^\prime\setminus p(\Supp R_{|S})$, where $\varepsilon_0$ is a positive real, and $\varphi_{\Lambda}$, respectively 
 $\varphi_{\Theta_{W^\prime}}$ are the local potentials of $\Lambda$, respectively $\Theta_{W^\prime}$;
\smallskip
\item {\rm (c)} The support of $\Lambda$ is contained in $\cup_j F_j$. }
\hfill\qed

\endclaim
\medskip

\noindent The point b) is the analogue of the klt property ii) in the result quoted above. Before explaining the main ideas involved in our arguments, we make here a few observations.

\claim Remarks|{\rm 

{
\itemindent 4mm

\item{1)} By the $L^2$ theory we can``convert'' the current $\Theta_{W^\prime}$ into an effective divisor
$\Delta_{W^\prime}$, as follows. For each $\varepsilon > 0$, there exists 
an effective divisor $\Delta_{W^\prime}$ on $W^\prime$ such that the part a) of the statement above becomes
the familiar relation
$$K_{W^\prime}+ \Delta_{W^\prime}\equiv \Lambda+ g^\star ({K_X+ \Delta+ \varepsilon A})_{|W^\prime};$$
and such that b) is still satisfied.
The precise procedure will be detailed at the end of our article.
As we have already mentioned, the minimal center $W$ is
known to be normal, provided that 
the pair $(X, \Delta)$ belongs to an appropriate class of singularities
(see e.g. the works of F. Ambro and O. Fujino in [1], [13]). Hence in this case
the divisor $\Lambda$ is contractible with respect to the map $g$, which means that by projection to $W$, the relation above is the same as the point $\rm i)$
of the original theorem of Kawamata
recalled above.  

\item{2)}  We can 
take $\Lambda= 0$ if the canonical bundle $K_{W^\prime}$ happens to be the inverse image of some bundle on $W$, as a by-product of our proof.

\item{3)} In the approach by J. Koll\'ar in [23], (the analogue of) the difference 
$\Theta_{W^\prime}- \Lambda$
is further analyzed: he exhibits a nef part, as well as a klt and canonical part
of it. Unfortunately, our methods as they stand by now do not allow to obtain such a decomposition
(it is of course not impossible that one could get it by a more refined version of theorem 0.1).\hfill\qed

}
}
\endclaim

\medskip
\noindent Next we will give a brief account of our proof of 0.3.
We first observe that we have
$$p^\star\big(g^\star(K_X+ \Delta)_{|W^\prime}-K_{W^\prime}\big)+ \Xi_{|S}\equiv K_{S/W^\prime}+ \Delta^\prime_{|S}+  R_{|S}.$$

Let $m_0$ be a positive and divisible enough integer; the main point of our strategy is to  show that the current $m_0\Theta_{S/W^\prime}$ associated to the 
$m_0$-Bergman kernel metric of the bundle 
$$m_0(K_{S/W^\prime}+ \Delta^\prime_{|S})$$
is at least as singular as the divisor type current $m_0[\Xi_{2|S}]$.
Very roughly, this is proved to be so in two steps. We first observe that thanks to the 
numerical relation above together with (3), we can ``move" the current $m_0\Theta_{S/W^\prime}$ into the class
$m_0(K_{S}+ \Delta^\prime_{|S})+ m_1\mu^\star (A)_{|S}$ (for some large enough positive integer $m_1$) in such a way that the only {\sl additional singularities} 
we get by this operation are concentrated along the support of $\Xi_{1|S}$. We approximate this new current with global sections of $km_0(K_{X^\prime}+ S+ \Delta^\prime)+ km_1\mu^\star (A)+ A^\prime_{|S}$, where
$k\gg 0$ and $A^\prime$ is an ample line bundle on $X^\prime$ (which does not depends on $k$), by using the {\sl regularization of closed positive currents}. 
The second step consists in showing that the approximating sections, twisted with
$km_0R_{|S}$ extend to $X^\prime$; here we use the 
{\sl invariance of plurigenera technique} originally due to Y.-T. Siu in [32], [33], in a slight variation of the version proved by L. Ein and M. Popa in [12]. 
This result is absolutely central for our proof, but nevertheless
we have decided to include its proof in the sequel of this article [7], since the arguments needed for it are quite different from the topics presented here.
The claim concerning the singularities acquired by $m_0\Theta_{S/W^\prime}$
follows from the 
Hartogs principle.

The existence of this particular representative of the 
Chern class of the right hand side in the formula above shows that the divisor
$K_{S/W^\prime}+ \Delta^\prime_{|S}+  R_{|S}-\Xi_{2|S}$ is
pseudoeffective. The current
$\Theta_{W^\prime}$ and respectively the divisor $\Lambda$ are obtained by taking the direct image of this difference.
We would like to stress on the fact that the explicit expression 
of the $m_0$-Bergman kernel metric in the theorem 0.1 is crucial for 
the proof of the integrability property b).\hfill \qed

\claim Acknowledgments|{\rm Part of the questions addressed in this article were
brought to our attention by L. Ein and R. Lazarsfeld; we are very grateful 
to them for this. We would also like
to thank F. Ambro, Y. Kawamata, J. Koll\'ar, M. Popa and E. Viehweg for interesting and stimulating
discussions about various subjects. Finally, we gratefully acknowledge
support from the Mittag-Leffler institute.

}
\endclaim

\bigskip
 
 \vskip 10pt

\section {\S A. Proof of the main result}
\vskip 10pt We start by
establishing some technical tools which will be needed later on in our arguments. In the paragraph A.1 we deal with a local 
$L^{2/m}$ extension result (proposition 0.2), whereas in A.2  we prove the theorem 0.1. 

\subsection {\S A.1 The $L^{2/m}$ version of the Ohsawa-Takegoshi theorem} 

\vskip 5pt 
In this section we  establish an extension theorem with $L^{2/m}$
bounds, 
analogous to the Ohsawa-Takegoshi result; as we have already mentioned, we will need it for the proof of the theorem 0.1. We recall that the setting is the following: let $\Omega\subset \bC^n$ be a ball of radius $r$ and let $h:\Omega\to \bC$ be a holomorphic function, such that 
$\sup_\Omega|h|\leq 1$; moreover, we assume that the gradient $\partial h$ of $h$ 
is nowhere zero on the set $V:=  (h=0)$. We denote by $\varphi$ a plurisubharmonic function, such that its restriction to $V$ is well-defined (i.e., $\varphi_{|V}\not \equiv -\infty$).
\medskip 
\proof $\hskip 1mm${\sl of 0.3}.
Let $f$ be a  holomorphic 
function
$f: V\to \bC$  with the property that 
$$\int_{V} |f|^{2/m}\exp (-\varphi){{d\lambda_V}\over {|\partial h|^2}}= 1. $$ 
We begin with some reductions. In the first place we can assume that the function $\varphi$
is smooth, and that the functions $h$ (respectively $f$) can be extended in a neighbourhood of $\Omega$ (of $V$ inside $V\cap \ol\Omega$, respectively). Once the result 
is established under these additional assumptions, the general case follows by approximations and standard normal families arguments.

We can then clearly find some holomorphic $F_1$ in $\Omega$ that extends $f$ and satifies
$$
\int_\Omega |F_1|^{2/m}\exp(-\varphi)d\lambda\leq A<\infty.
$$
We then apply the Ohsawa-Takegoshi theorem with weight
$$
\varphi_1=\varphi +(1-1/m)\log|F_1|^2
$$
and obtain a new extension $F_2$ of $f$ satisfying
$$
\int {{|F_2|^2}\over{|F_1|^{2-2/m}}}\exp-\varphi\, d\lambda\leq
C_0\int {{|f|^2}\over{|F_1|^{2-2/m}}}\exp-\varphi \,{{d\lambda_V}\over{|\partial h|^2}}
=C_0.
$$
H\"older's inequality gives that
$$
\int_\Omega |F_2|^{2/m}\exp(-\varphi)\, d\lambda \int_\Omega {{|F_2|^{2/m}}\over{|F_1|^{(2-2/m)/m}}}|F_1|^{(2-2/m)/m}\exp(-\varphi) \,d\lambda\leq 
$$
$$
\leq(\int_\Omega {{|F_2|^2}\over{|F_1|^{2-2/m}}}\exp(-\varphi) d\lambda)^{1/m}
(\int_\Omega|F_1|^{2/m}\exp(-\varphi)d\lambda)^{m/(m-1)}
$$
which is smaller than
$$
C_0^{1/m} A^{m/(m-1)}= A(C_0/A)^{1/m}=:A_1.
$$
If $A>C_0$, then $A_1<A$. We can then repeat the same argument with $F_1$ replaced by $F_2$, etc, and get a decreasing sequence of constants $A_k$, such that 
$$A_{k+1}= A_k(C_0/A_k)^{1/m}$$
for $k\geq 0$. It is easy to see that $A_k$ tends to $C_0$. Indeed, if $A_k >r C_0$ for some $r>1$, then $A_k$ would tend to zero by the relation above. This completes the proof.
\hfill \qed

\vskip 10pt
\noindent For further use, we reformulate the above result in terms of differential forms. One can find the details of the ``conversion" of 0.2 into $0.2^\prime$ e.g. in the notes [5].

\claim $0.2^\prime$ Lemma|Given $u$  a holomorphic $m$-canonical form on $V$ such that 
$$\int_{V} |u|^{2\over m}e^{-\varphi}<\infty, $$
there exist a holomorphic $m$-canonical form $U$ on $\Omega$ such that:
\item{(i)} $U= u\wedge (dh)^{\otimes m}$ on $V$ i.e. the form $U$ is an extension of $u$.
\smallskip
\item{(ii)} The next $L^{2/m}$ bound holds 
$$\int_{\Omega}  |U|^{2\over m}e^{-\varphi} \leq C_0\int_{V}  |u|^{2\over m}e^{-\varphi}, $$where $C_0$ is the constant appearing in the Ohsawa-Takegoshi theorem.
\hfill\qed
\endclaim
\bigskip

\subsection{\S A.2 Positivity properties of the $m$-Bergman kernel metric}

\smallskip

For the convenience of the reader, we first recall some facts concerning the
variation of the fiberwise 
$m$--Bergman metric (also called NS metric in [6], [30], or just a pseudo-norm in [17]).

The general set-up is the following: let $p: X\to Y$ be a surjective projective map between
two non-singular manifolds and consider  a line bundle $L\to X$, endowed with a 
metric $h_L$, such that $\sqrt {-1}\Theta_h(L)\geq 0$. To start with, we will assume 
that {\sl $h_L$ is a genuine metric (i.e. non-singular) and that the map $p$ is a smooth fibration. }

With this data, for each positive integer $m$ we can construct a 
metric $\displaystyle h^{(m)}_{X/Y}$ on the twisted relative bundle 
$\displaystyle mK_{X/Y}+ L$ as follows: let us take a vector $\xi$ in the fiber  
$-(mK_{X/Y}+ L)_x$; then we define its norm 
$$
\Vert \xi\Vert ^2: =\sup |\wt \sigma(x)\cdot \xi|^2
$$
the "$\sup$`` being taken over all sections $\sigma$ to $\displaystyle mK_{X_y}+L$ such that
$$
\Vert \sigma\Vert_{m, z}^{2/m}:= \int_{X_y} |\sigma|^{2\over m}e^{-{1\over m}\varphi_L}\leq 1;
$$
in the above notation $\varphi_L$ denotes as usual the metric of $L$ and $p(x)= y$; also, we have used the 
notation $\wt \sigma$ to denote the section of the bundle $\displaystyle mK_{X/Y}+ L_{|X_y}$ corresponding to
$\sigma$ via the standard identification (see e.g. [6], section 1). 
\vskip 5pt
We  give now the expression of the local weights of the metric 
$h^{(m)}_{X/Y}$; we denote by $(z^j)$ respectively $(t^i)$ some local coordinates centered at $x$, respectively $y$. Then we have
$$\exp \big(\varphi^{(m)}_{X/Y}(x)\big)= \sup _{\Vert u\Vert _{m, y}\leq 1}\vert u^\prime(x)|^2$$
where $u\in H^0\big(X_y, mK_{X_y}+L\big)$ is a global section.

Finally, we denote by $\wt u:= u\wedge (p^*dt)^{\otimes m}$ and the above $u^\prime$ is obtained as  $\wt u= u^\prime (dz)^{\otimes m}$.
\medskip
 In the article [6], among other things we have established the fact that {\sl the metric $h^{(m)}_{X/Y}$ above has semi-positive curvature current, 
or it is identically $+\infty$}.
The latter situation occurs precisely when there are no global $L$--twisted $m$-canonical forms on the
fibers. 
\medskip
\noindent In fact, it turns out that the above construction has a meaning even if
the metric $h_L$ we start with is allowed to be singular (but we still assume that the map $p$ is a non-singular fibration). We remark that in this case some fibers $X_y$ may be contained in the unbounded locus of $\displaystyle h_L$, i.e.
$\displaystyle h_{L|X_y}\equiv \infty$, but for such $y\in Y$ we adopt the convention that the 
metric $h^{(m)}_{X/Y}$ is identically $+\infty$ as well. As for the fibers in the complement of this set, 
the family of sections we consider in order to define the metric consists in twisted pluricanonical forms whose $m^{th}$ root is $L^2$. 
\medskip
\noindent We recall now the the result we have proved in the section 3 of [6] (see also [30] and the references therein).
\claim A.2.1 Theorem {([6])}|Let $p: X_0\to Y_0$ be a proper projective non-singular fibration, and let 
$L\to X_0$ be a line bundle endowed with a metric $h_L$ such that:

{\itemindent 4mm
\smallskip
\item {\rm 1)} The curvature current of the bundle $(L, h_L)$ is positive, i.e. 
$\displaystyle \sqrt {-1}\Theta_{h_L}(L)\geq 0$;

\item {\rm 2)} For each $y\in Y_0$, all the sections of the bundle 
$mK_{X_y}+ L$
extend near $y$;

\item {\rm 3)} There exist $z\in Y_0$ and a section 
$\displaystyle u\in H^0\big(X_z, mK_{X_z}+ L)\big)$
such that 
$$\int_{X_y} |u|^{2\over m}e^{-{1\over m}\varphi_L}< \infty.$$

}

\noindent Then the above fiberwise $m$--Bergman kernel metric $\displaystyle h^{(m)}_{X_0/Y_0}$ has semi-positive
curvature current.\hfill\qed
\endclaim
\medskip

\claim A.2.2 Remark|{\rm We take this opportunity to point out two {\sl imprecisions} 
in the formulation of the corollary 4.2 from our previous work [6]. 
In the first place, it must be of course assumed that the metric $\displaystyle h^{(m)}_{X_0/Y_0}$ (called NS-metric in [6]) is not identically $+\infty$. 
Secondly, the assumption that $p_\star(mK_{X/Y}+ L)$ be locally free should be replaced by 
the assumption 2) of the preceding statement. Immediately after corollary 4.2 in [6] it is 
stated that these two conditions are equivalent, but this is not the case. In fact, the property of being locally free concerns the direct image of $p_\star(mK_{X/Y}+ L)$ as a sheaf: the stalks of this sheaf consists of sections over the fibers of $p$ that {\sl do extend} to neighboring
fibers, and local freeness means that this space of extendable sections has everywhere the same rank. We remark that the direct image $p_\star(mK_{X/Y}+ L)$ is torsion free, hence locally free if the dimension of the base is equal to one. For the proof of Corollary 4.2 we need the stronger property 2) above, saying that all sections extend locally to neighboring fibers.

This distinction between local freeness and the condition 2) is precisely the heart of the question of invariance of plurigenera for projective families over one dimensional base, where the local freeness is automatic, and 2) is exactly what is to be proved. We will continue this discussion at the end of the present subsection, after the proof of the corollaries of the Theorem 0.1. \hfill\qed

}
\endclaim

\medskip
\proof \hskip 1.6mm (of 0.1). By the usual semi-continuity arguments, there exist a Zariski open set $Y_0\subset Y$ such that the condition (2) above is fulfilled; by restricting further the set $Y_0$ we can assume that if we denote by $X_0$ the inverse image 
of the set $Y_0$ via the map $p$, then $p:X_0\to Y_0$ is a non-singular fibration.

Then we can use the above result and infer that our metric $h^{(m)}_0$ is explicitly given over the fibers $X_y$, as soon as $y\in Y_0$ and the restriction of $h_L$ to $X_y$ is well defined.
We are going to show now that this metric admits an extension to the whole manifold $X$. The method of proof is borrowed from our previous work: we show that the local weights of the metric
$h^{(m)}_0$ are locally bounded near $X\setminus X_0$ and then standard results in pluripotential theory imply the result.

Let us pick a point $x\in X$, such that $y:= p(x)\in Y_0$; assume that the restriction of $h_L$ to $X_y$ is well defined. Let $u$ be some
global section of the bundle $\displaystyle mK_{X_y}+ L$. Locally near $x$ we consider a coordinates ball $\Omega$ and let $\Omega_y$ be its intersection with the fiber $X_y$. 
Thus on $\Omega_y$ our section $u$ is just the $m$'th tensor power of some $(n, 0)$ form, and let 
$\wt u:= u\wedge p^{\star}(dt)^{\otimes m}$. With respect to the local coordinates $(z^{j})$ 
on $\Omega$ we have $\wt u= u^\prime dz^{\otimes m}$ and as explained before, the local weight 
of the metric $h^{(m)}_0$ near $x$ is given by the supremum of $|u^\prime|$ when $u$ is
normalized by the condition
$$\int_{X_y} |u|^{2\over m}e^{-{1\over m}\varphi}\leq 1.$$

But now we just invoke the lemma $0.2^\prime$, in order to obtain some $m$--form of maximal degree 
$\wt U$ on $\Omega$ such that its restriction on $\Omega_y$ is $\wt u$, and such that the next inequality hold
$$\int_\Omega |\wt U|^{2\over m}e^{-{1\over m}\varphi}\leq C_0 
\int_{\Omega_y}|u|^{2\over m}e^{-{1\over m}\varphi}.$$

The mean value inequality applied to the relation above show that we have  
$$|u^\prime(x)|^{2/m}\leq C\leqno (4)$$ 
where moreover the bound $C$ above does not depend at all on the geometry of the fiber $X_y$, but on the ambient manifold $X$. In particular, the metric 
$h^{(m)}_0$ admits an extension to the whole manifold $X$ and its curvature current is semipositive, so our theorem 0.1 is proved. \hfill\qed
\bigskip
\noindent {\bf A.2.3 Remark.} The constant $``C "$ above
only depends on the sup of the $m^{th}$--root of the metric of the bundle $L$ and the geometry of the ambient manifold $X$ (and {\sl not at all} on the geometry of the fiber $X_y$). As a consequence we see that if we apply the above arguments to a sequence
$pK_{X/Y}+ L_p$, then the constant in question will be uniformly bounded 
with respect to $p$, provided that $\displaystyle {1\over
  p}\varphi_{L_p}$ is bounded from above. Hence we obtain the corollaries below (see also [11] and [30] for related statements).

\claim A.2.4 Corollary|Let $p: X\to Y$ be a projective surjective map, and let 
$L\to X$ be a line bundle endowed with a metric $h_L$ such that:

{\itemindent 4mm
\smallskip
\item {\rm 1)} The curvature current of the bundle $(L, h_L)$ is positive, i.e. 
$\displaystyle \sqrt {-1}\Theta_{h_L}(L)\geq 0$;

\item {\rm 2)} There exists a very general point $z\in Y$ and a section 
$\displaystyle u\in H^0\big(X_z, m_0K_{X/Y}+ L_{|X_z}\big)$
such that 
$$\displaystyle \int_{X_z} |u|^{2\over m_0}e^{-{1\over m_0}\varphi_L}< \infty.$$

}

\noindent Then the bundle  $mK_{X/Y}+ L$ admits a metric $\varphi^{(\infty)}_{X/Y}$ with 
positive curvature current, whose restriction to any very general fiber $X_w$ of $p$ 
has the following (minimality) property: for any non-zero section $\displaystyle v\in H^0\big(X_w, kmK_{X_z}+ kL\big)$ we have
$$|v|^{2\over km}e^{-\varphi^{(\infty)}_{X/Y}}\leq
\int_{X_w} |v|^{2\over km}e^{-{1\over m}\varphi_L}$$
up to the identification of $K_{X/Y|X_w}$ with $K_{X_w}$.
\hfill\qed

\endclaim
  
\medskip 

\noindent In the statement A.2.4 we call a point $z\in Y$ {\sl very general} if any holomorphic section of the bundle
$$km_0K_{X/Y}+ kL_{|X_z}$$
extends locally near $z$, for any positive integer $k$.
By the usual semi-continuity arguments we infer that 
the set of very general points is the complement of a countable union
of Zariski closed sets of codimension at least one in $Y$. \hfill\qed
\smallskip

\noindent The  proof of the above result is an immediate consequence of the arguments of 0.1. Indeed, we consider the $km$--Bergman metric on the bundle
$$k\big(mK_{X/Y}+ L\big)$$
and we know that the weights of its $1/k^{\rm th}$ root are {\sl uniformly} bounded, independently of $k$ (see the remark above). We obtain the metric $\varphi^{(\infty)}_{X/Y}$ by the usual upper envelope construction, namely
$$\varphi^{(\infty)}_{X/Y}: = {\sup} _{k\geq 0}^\star \Big( {1\over k}\varphi^{(km)}_{X/Y}\Big).$$
\hfill\qed

\medskip
\noindent The corollary A.2.4 admits a metric version which we discuss next.

\claim A.2.5 Corollary|Let $p: X\to Y$ be a projective surjection, where $X$ and 
$Y$ are non-singular manifolds, and let $L\to X$ be a line bundle, 
endowed with a positively curved metric $h_L$. 
We assume that there exists a very general point $w\in Y$ and a metric $h_{w}$
of the bundle $m_0K_{X_w}+ L$ with positive curvature and such that 
$$e^{\varphi_w-\varphi_L}\in L^{{{1+\varepsilon}\over {m_0}}}\leqno (5)$$
holds locally at each point of $X_w$, where $\varepsilon> 0$ is a positive real.. 

\noindent Then the bundle $m_0K_{X/Y}+ L$ admits a 
positively curved metric whose restriction to $X_w$ is less singular than the metric $h_{w}$.
\endclaim
\noindent Let $A= 2nH$ where $n$ is the dimension of $X$ and $H$ is a
very ample bundle on $X$.
In the theorem above, a point 
$y\in Y$ is called generic if it is not a critical value of $p$ and if 
any section of the bundle $km_0K_{X/Y}+ kL+A_{|X_w}$ extends 
over the nearby fibers, for all $k$.

We remark that the corollary A.2.5 is more {\sl coherent} than A.2.4, 
in the sense that we start with a metric on
a very general fiber $X_w$ of $p$, and we produce a metric defined on $X$.
Even if we give ourselves a section on a very general fiber, the object 
we are able to produce (via A.2.4) will be in general a metric, which is not necessarily algebraic, as one would hope or guess. 

We first write $e^{-\varphi_w}$ as limit of a sequence of algebraic metrics $e^{-\varphi_w^{(k)}}$ of
$km_0K_{X/Y}+ kL+ A_{|X_w}$,
and then we apply 0.1 to construct a sequence of global metrics 
on 
$$km_0K_{X/Y}+ kL+ A,$$
whose restriction to $X_w$ is comparable with $e^{-\varphi_w^{(k)}}$.
The metric we seek is obtained by a limit process. In order to complete this program 
we have to control several constants which are involved in the arguments, 
and this is possible 
thanks to the version of the Ohsawa-Takegoshi theorem 0.2. 
It would be very interesting to give a more direct, {\sl sections-less} proof of A.2.5. \hfill\qed
\medskip

\proof \hskip 1.6mm (of A.2.5).  The choice of the bundle $A\to X$ as above is explained by the following fact we borrow from the approximation theorem in [9]. 

Associated to the metric $\varphi_w$ on the bundle $m_0K_{X/Y}+ L_{|X_w}$ given by hypothesis
we consider the following space of sections
$$\cV_k\subset H^0(X_w, km_0K_{X_w}+ kL+ A)\leqno (6)$$ defined by
$u\in \cV_k$ if and only if
$$\Vert u\Vert_w^2:= \int_{X_w}|u|^2\exp(-k\varphi_w-\varphi_A)dV_\omega< \infty\leqno (7)$$
We denote by $\varphi_w^{(k)}$ the metric on the bundle $km_0K_{X_w}+ kL+ A$ induced by an orthonormal basis $\big(u^{(k)}_{j}\big)$ of $\cV_k$ endowed with the scalar product corresponding to (2). Then it is proved in [9] that there exists a constant $C_1$ independent of $k$ such that we have 
$$k\varphi_w+ \varphi_A \leq \varphi_w^{(k)}+ C_1.\leqno (8)$$ 
We remark that the ample bundle $A$ is independent of the particular point 
$w\in Y$ as well as on the metric $e^{-\varphi_w}$. This fact is very important for the 
following arguments.

\noindent 
The H\"older inequality shows that the 
$km_0^{\rm th}$ root of the sections $v\in \cV_k$ are integrable with respect to the metric 
$\displaystyle {{1}\over {m_0}}\varphi_L$, as soon as $k$ is large enough. Indeed we have:

$$\int_{X_w}|v|^{2\over km_0}e^{-{k\varphi_L+ \varphi_A\over km_0}}\leq
c_k\Big(\int_{X_w}|v|^2e^{-k\varphi_w- \varphi_A}dV_\omega\Big)^{{{1}\over {km_0}}}\leqno (9)$$
where we use the following notation
$$c_k^{{{km_0}\over {km_0-1}}}:= \int_{X_w}e^{{{k}\over {km_0-1}}(\varphi_w
-\varphi_L)}dV_\omega< \infty.$$
We remark that the last integral is indeed convergent, by the integrability hypothesis
in (5), provided that $k\gg 0$. In conclusion, we get
$$\int_{X_w}|v|^{2\over km_0}e^{-{k\varphi_L+ \varphi_A\over km_0}}\leq C_2
\leqno (10)$$
for every section $v:= u^{(k)}_{j}\in \cV_k$. By definition, the expression of the constant in the relation (10) is the following
$$C_2:= \sup_{k\gg 0} c_k< \infty $$ therefore it is independent of $k$.
 \smallskip
By theorem 0.1, for each $k$ large enough 
we can construct the $km_0$-Bergman metric $h^{(km_0)}_{X/Y}$ on the bundle
$km_0K_{X/Y}+ kL+ A$ with positive curvature current. Thanks to the relation (10), its restriction to 
$X_w$ is less singular than the metric $\varphi^{(k)}_w$, in the following precise way. 
$$\varphi^{(km_0)}_{X/Y|X_w}\geq \varphi^{(k)}_w- km_0\log C_2;\leqno (11)$$
(we identify $K_{X/Y}$ with $K_{X_w}$, which is harmless in this context, since the point
$w$ is ``far'' from the critical loci of $p$)
\noindent The inequality (11) combined with (8) shows furthermore that we have
$${1\over k}\varphi^{(km_0)}_{X/Y|X_w}\geq \varphi_w+ C_4\leqno (12)$$
where $C_4:= -m_0\log C_2$.

On the other hand, by remark A.2.3 we have an a-priori upper bound 
$${1\over k}\varphi^{(km_0)}_{X/Y}\leq C_5\leqno (13)$$
where $C_5$ is {\sl uniform with respect to $k$}. 

The metric we seek is obtained by the usual upper envelope construction, namely
$$\varphi^{(\infty)}_{X/Y}: = \lim {\sup} _{k\geq 1}^\star \Big({1\over k}\varphi^{(km_0)}_{X/Y}\Big).\leqno (14)$$
The local weights $\big(\varphi^{(\infty)}_{X/Y}\big)$ glue together to give a metric
for the bundle $m_0K_{X/Y}+ L$, since the auxiliary bundle $A$ used in the approximation process
is removed by the normalization factor $1/k$.

Moreover, the relation (12) shows that we have 
$$\varphi^{(\infty)}_{X/Y|X_w}\geq \varphi_w+ \cO(1); \leqno (15)$$
hence, the metric $\varphi^{(\infty)}_{X/Y}$ restricted to the fiber $X_w$ is 
less singular than the metric $e^{-\varphi_w}$ and the corollary is proved.
\hfill\qed

\medskip
\claim A.2.6 Remark|{\rm It is certainly worthwhile to understand the behavior of the metric constructed in Theorem 0.1 over the {\sl exceptional fibers}, where it is just defined as the (unique) extension of $h^{(m)}_{X_0/Y_0}$.

Let us here look at the fibers over a point $a$ where the condition 2) in 0.1 is not known to be satisfied, but such that $a$ is still a regular value of $p$. On the bundle
$\displaystyle m_0K_{X_a}+ L$ we can consider (at least) two natural extremal metrics:
the one induced by {\sl all} its global sections satisfying the $L^{2/m}$ integrability 
condition with respect to $h_L$ which we denote by $h_1$, and metric corresponding to the subspace of sections which extend locally near $a$, denoted by $h_2$.
It is clear that we have $h_2\geq h_1$, as one can see by the comparison between the corresponding $m$-Bergman kernels. 

\noindent Our claim is that the metric $h^{(m)}_{X_0/Y_0}$ extends over $X_a$ to a metric
which is {\sl less singular than} $h_2$, at least if the singularities of the metric 
$\varphi_{L}$ are mild enough. 

Indeed, as a consequence of the theorem 0.1 we have
$$\varphi_{X/Y}^{(m)}(x)= \lim\sup_{x^\prime\to x}\varphi_{X_0/Y_0}^{(m)}(x^\prime)$$
where $x^\prime\in X_0$ in the limit above. In order to establish the result we claim, let us consider a holomorphic section $u\in H^0(X_a, m_0K_{X_a}+ L)$ which computes the local weight of $h_2$ at $x$, i.e.
$\Vert u\Vert_{m, a}= 1$ and 
$$\sup_{\Vert v\Vert_{m, a}= 1}|v^\prime(x)|^2= |u^\prime(x)|^2$$
in the notations/conventions at the beginning of the section A.2. If $x^\prime\in X_0$
is close enough to $x$ and $y^\prime:= p(x^\prime)$, then we certainly have
$$\sup_{\Vert v\Vert_{m, y^\prime}=1}|v^\prime(x^\prime)|^2\geq 
 {{|u^\prime(x^\prime)|^2}\over {\Vert u\Vert_{m, y^\prime}^2}}$$
(here we use the fact that the section $u$ extends over the fibers near $X_a$). 
Hence in order to establish the claim above, all that we need is the relation
$$\lim\inf_{y^\prime\to a}\Vert u\Vert_{m, y^\prime}\leq \Vert u\Vert_{m, a}.$$
This last inequality clearly holds at least if $\varphi_L$ is continuous, and we actually believe
that {\sl it always holds}. \hfill\qed 
}\endclaim


\bigskip
\section {\S B. Adjunction type results}

\medskip 
\noindent The manner in which the technics developed in the first part of our article are used in order to prove 0.3 is explained in the subsection B.1 and B.2, together with comments about more general statements.

We use the notations in the introduction.
A first step in the direction of our proof is to write the relation $(\dagger)$ as follows:

$$p^*\big(g^\star (K_X+ \Delta_{|W})\big)+ \Xi_{|S}\equiv 
K_{X^{\prime}}+ S+ \Delta^{\prime}+  R_{|S}\leqno (16)$$
We subtract the inverse image of the canonical bundle of $W^\prime$ and we get
$$p^*\big(g^\star (K_X+ \Delta_{|W})-K_{W^\prime}\big)+ \Xi_{|S}\equiv 
K_{S/W^\prime}+ \Delta^{\prime}+  R_{|S}\leqno (17)$$	 
by restriction to $S$.

Let $m_0$ be a positive integer, such that the multiples $m_0(K_X+ \Delta)$, $m_0\Delta^{\prime}$ and $m_0R$ are
integral. 
Given the formula (17) above and the statement we want to obtain, it is clear that our first task will be to analyze the positivity properties of the bundle
$$m_0(K_{S/{W^\prime}}+ \Delta^\prime_{|S});\leqno (18)$$
this will be done in the next paragraph.
\bigskip

\subsection {\S B.1 The Bergman metric and its singularities}

\medskip 
\noindent We state here the following direct consequence of the theorem 0.1.

\claim B.1.1 Corollary|The bundle $m_0(K_{S/{W^\prime}}+ \Delta^\prime_{|S})$ is pseudoeffective.
Moreover, it admits a positively curved metric $h_{S/{W^\prime}}$,
whose restriction to the generic fiber
$S_{w^\prime}$ of $p$ is induced by the space of sections $H^0\big(S_{w^\prime}, {m_0(K_{S/{W^\prime}}+ \Delta^\prime}_{|S_{w^\prime}})\big
)$ as explained in the paragraph A.2.
\endclaim

\proof. Let $L$ be a line bundle whose Chern class contains the 
Weil divisor $m_0\Delta^\prime_{|S}$, and such that the numerical relation (17) becomes 
linear when we replace $\Delta^\prime_{|S}$ with $1/m_0L$.
The bundle $L$ admits a singular metric 
$h_L$ whose curvature form is equal to $m_0[\Delta^\prime_{|S}]$.
The multiplier ideal sheaf corresponding to the $m_0^{\rm th}$ root of this metric is  
equal to the structural sheaf, since the coefficients $a^j$ belong to the interval $[0, 1[$
by definition of $\Delta^\prime$. Therefore, given a generic
point $w^\prime\in W^\prime$, we infer that any section of the bundle 
${m_0K_{S/{W^\prime}}+ L}_{|S_{w^\prime}}$
will automatically satisfy the condition 2 of 0.1.
\smallskip
Furthermore, we remark that we have 
$$H^0\big(S_{w^\prime}, {m_0K_{S/{W^\prime}}+ L}_{|S_{w^\prime}}\big
)\neq 0$$
since the support of the divisor $R_{|S}$ {\sl does not} intersects $S_{w^\prime}$, and since the bundle on the left hand side of $(17)$ is effective
when restricted to the generic 
fiber of $p$. Thus we obtain a well-defined $m_0$--Bergman kernel metric $h_{S/{W^\prime}}$
on the bundle $m_0K_{S/{W^\prime}}+ L$, which proves the corollary. \hfill\qed
\medskip
\noindent During the rest of this section, we use the symbol $\displaystyle\Theta_{S/{W^\prime}}$ to denote $1/m_0$ times the curvature current associated to $h_{S/{W^\prime}}$. We remark that $\displaystyle\Theta_{S/{W^\prime}}$ is closed, positive, and 
it belongs to the cohomology class  $\{K_{S/{W^\prime}}+ \Delta^\prime_{|S}\}$.Å\medskip  
\noindent  
As a consequence of B.1.1, we note that by (17) and $(\dagger \dagger)$ we have the next identity
$$p^*\big(g^\star (K_X+ \Delta)_{|W^\prime}-K_{W^\prime}\big)+ [\Xi_{1|S}]+ [\Xi_{2|S}]\equiv 
\Theta_{S/{W^\prime}}+  [R_{|S}].\leqno (19)$$
The existence of the current $\Theta_{S/{W^\prime}}$ is certainly
excellent news in view of what is to be proved in 0.3; however, 
in order to conclude, we have to show that this current is {\sl at least} as singular as the current $[\Xi_{2|S}]$.
Once this is done, the objects we seek $\big(\Theta_{W^\prime}, \Lambda\big)$ 
will be obtained by direct image.
\medskip 

\noindent The claim concerning the singularities of the current $\Theta_{S/{W^\prime}}$ is a consequence of the 
following statement.

\claim B.1.2 Theorem {\rm ([7])}|Let $T$ be any  closed positive $(1,1)$-current, for which there exists a line bundle $F$ on 
$W^\prime$ with the property that
$T\equiv p^\star (F)+ m_0(K_{S/{W^\prime}}+ \Delta^\prime_{|S})$. Then we have 
$$T\geq [\Xi_{2|S}]$$ 
in the sense of currents on $S$.
\hfill\qed

\endclaim
 
 \smallskip 
 \noindent Despite the fact that this result is absolutely central for what will follow, 
we have decided to include its proof in the sequel of this article [7], since the arguments needed for it are quite different from the topics presented here. \hfill\qed
 
\medskip
\noindent We now continue the proof of 0.3 by admitting B.1.2.
 
 The relation (19) can be re-written as follows
 $$p^*\big(g^\star (K_X+ \Delta_{|W})-K_{W^\prime}\big)+ [\Xi_{1|S}]\equiv 
\Theta_{S/{W^\prime}}- [\Xi_{2|S}]+  [R_{|S}],\leqno (20)$$
 where the right hand side term is a {\sl closed positive current}. The positivity assertion is a consequence of the theorem B.1.2, where $T:= \Theta_{S/{W^\prime}}$
and $F$ is the trivial bundle.
 Next, we define
$$\Lambda:= \sum_jq^jF_j\leqno(21)$$
to be the smallest effective $\bQ$-divisor such that 
$$p^\star(\Lambda)\geq \Xi_{1|S}.\leqno(22)$$
Such an object indeed exists, by the definition of the decomposition
$\Xi= \Xi_{1}+ \Xi_2$ in the introduction. We denote by $D:= p^\star(\Lambda)- \Xi_{1|S}$; it is an effective $\bQ$-divisor, and then the relation (20) reads as follows
$$p^*\big(g^\star (K_X+ \Delta_{|W})-K_{W^\prime}+ \Lambda\big)\equiv 
\Theta_{S/{W^\prime}}- [\Xi_{2|S}]+  [R_{|S}]+ [D].\leqno (23)$$
Again, we notice that the right hand side term of the preceding relation 
is a closed positive current. Since $\Theta_{S/{W^\prime}}- [\Xi_{2|S}]+  [R_{|S}]+ [D]$ 
belongs to the $p$-inverse image of a $\bQ$-bundle on $W^\prime$, we infer the 
existence of a closed positive current $\Theta_{W^\prime}$ on  $W^\prime$, such that 
$$g^\star (K_X+ \Delta_{|W})-K_{W^\prime}+ \Lambda\equiv 
\Theta_{W^\prime}\leqno (24)$$
which is precisely the relation we are looking for. \hfill\qed
\smallskip
In conclusion, the {\sl qualitative} part of 0.3 is proved; we turn in the next subsection to
the {\sl quantitative} part of it, namely the integrability statement b).

 \bigskip

\subsection {\S B.2 The critical exponent of $\Theta_{W^\prime}- \Lambda$}

\noindent  Our main goal in this paragraph is to prove the point b) of 0.3.
To this end,
we will first provide a closer analysis of the structure of the curvature current
$\Theta_{S/{W^\prime}}$ when restricted to the generic fiber of $p$. Our first statement
can be seen as a generalization of the result in [19], [23] corresponding to the case
$m_0= 1$.

\claim B.2.1 Theorem|Let $w\in W^\prime$ be a generic point and let 
$m\geq 1$ be a positive integer.
Then the divisor $\displaystyle mm_0\Xi_{|S_w}$ is contained in the zero set of any section
$u$ of the 
bundle 
$$m(m_0K_{S_w}+ L_{|S_w}).$$
Moreover we have $h^0\big(S_w, m(m_0K_{S_w}+ L_{|S_w})\big)= 1$.
\endclaim

\smallskip
\noindent {\bf Remark.} The bundle $L$ above is the one chosen in B.1.1.
The idea of the following proof is quite simple: we first show that the sections above twisted with inverse image of the section of some large enough ample line bundle on $W^\prime$ extends to $S$. Then we apply the theorem 
B.1.2 above.
\hfill \qed

\medskip
\noindent \proof \hskip 1.5mm (of B.2.1).  As in the proof of the corollary B.1.1, theorem 0.1 provides us with a metric $h^{(m)}_{S/W^\prime}$ on the bundle 
$\displaystyle m(m_0K_{S/W^\prime}+  L)$ with positive curvature current. 
If $w\in W^\prime$ is generic, then the restriction of the metric to the fiber $S_w$ has the same singularities
as  the algebraic metric given by the sections of the bundle 
$$\displaystyle m(m_0K_{S_w}+  L_{|S_w}).$$

\noindent Let $A^\prime\to W^\prime$ be a positive enough line bundle; we want to apply the Ohsawa-Takegoshi theorem to the bundle 
$m(m_0K_{S/W}+  L)+ p^\star A^\prime$ in order to extend $u\otimes p^\star(\sigma_A)$ to $S$, where $\sigma_A$ is a generic section of $A^\prime$.

To this end, we write it as an adjoint bundle as follows
$$m(m_0K_{S/W^\prime}+ L)+p^\star A^\prime= 
K_S+ \Big({mm_0-1\over m_0}\Big)\big(m_0K_{S/W^\prime}+ L\big)+{1\over m_0}L+ 
p^\star (A^\prime-K_{W^\prime})$$
The second bundle in the above decomposition is endowed with the appropriate multiple of the metric $h^{(m_0)}_{S/W^\prime}$. If we denote by 
$$F:= \Big({mm_0-1\over m_0}\Big)\big(m_0K_{S/W^\prime}+ L\big)+{1\over m_0}L+ p^\star 
(A^\prime-K_W^\prime)$$
then we have
$$u\otimes p^\star(\sigma_A) \in H^0\big(S_w, {K_S+ F}_{|S_w}\big)$$
and let $h_F$ be the metric on $F$ induced by the $h^{(m_0)}_{S/W^\prime}$ together with 
the natural metric on $\Delta_{|S}^\prime$, and with the smooth, positively curved metric
$h_A$ and on $p^\star(A^\prime-K_{W\prime})$.
Then we claim that the curvature assumptions and the $L^2$ conditions required by the 
Ohsawa-Takegoshi extension theorem are satisfied. 

To verify this claim, let $B\to W^\prime$ be a very ample line bundle, such that there exists a family of sections $\rho_j \in H^0(W^\prime, B)$ with the property that
$w= \cap_j(\rho_j= 0)$. The curvature conditions to be fulfilled are:

{\itemindent 4mm

\item {i)} $\displaystyle {{\sqrt {-1}}\over {2\pi}}\Theta_{h_F}(F)+{{\sqrt {-1}}\over {2\pi}}\ddbar \log(\sum |\rho_j\circ \mu|^2)\geq 0$
where we measure the norm of the sections $\rho_j$ above by a positively curved metric of $B$;
\smallskip
\item {ii)} $\displaystyle {{\sqrt {-1}}\over {2\pi}}\Theta_{h_F}(F)+ {{\sqrt {-1}}\over {2\pi}}\ddbar \log(\sum |\rho_j\circ \mu|^2)\geq \delta \mu^\star\big({{\sqrt {-1}}\over {2\pi}}\Theta (B)\big)$ where $\delta$ is a positive real.

}
\smallskip
\noindent We see that both conditions will be satisfied if the curvature of $F$
is greater than say $2\mu^\star\Theta (B)$; this last condition can be easily satisfied if $A$ is positive enough.

Concerning the integrability of the section, we remark that we have
$$\int_{S_w}|u|^2e^{-(m-{1\over m_0})\varphi^{(m_0)}_{S/W^\prime}-\varphi_{\Delta^\prime}-\varphi_{A}
}\leq C\int_{S_w}e^{-\varphi_{\Delta^\prime}}<\infty.$$
The first inequality is satisfied because $|u|^{2/m}$ is pointwise bounded with respect to 
$h^{(m_0)}_{S/W}$ and the second one holds because the multiplier ideal 
sheaf of $\displaystyle\Delta_{|S}^\prime$ restricted to $S_w$
is trivial, provided that $w$ is generic (which is assumed to be the case).

Then the Ohsawa-Takegoshi theorem shows that $u$ admits an extension 
$U$ as section of the bundle $m(m_0K_{S/W}+L)+ p^\star A^\prime$. 
But if this is so, we just invoke 
the theorem B.1.2 and infer that the divisor of zeroes of the section $u$ is greater than 
$mm_0\Xi_{2|S_{w^\prime}}$. 

In conclusion, we have just proved that the bundle $G:= m(m_0K_{S_w}+ L_{|S_w})$ 
has the next property: it has a non-zero holomorphic section $v$ such that for any other section $u$ 
of $G$, the quotient $\displaystyle u\over v$ is holomorphic. Since$\displaystyle u\over v$ is a section of the trivial bundle, it is equal to a constant
and thus theorem B.2.1 is completely proved. \hfill\qed

\vskip 5pt

\noindent If $D$ is an effective Weil divisor, we denote by $u_D$ the corresponding  section of the bundle $\cO(D)$. We consider a coordinate set $V\subset W^\prime$; 
there exists a {\sl meromorphic} section $u_{V}$ of the bundle
$${m_0K_{S/W^\prime}+ L}_{|p^{-1}(V)}\leqno (25)$$
defined by the relation 
$$\displaystyle u_{V}:= {{u_{m_0\Xi_{|S}}}\over {u_{m_0R_{|S}}}}.\leqno(26)$$
Here we implicitly use the fact that the bundle 
$${g^\star (K_X+ \Delta)_{|W^\prime}-K_{W^\prime}}_{|V}$$
is trivial, in order to identify the quotient in (26) with a section of the bundle (25). 

We remark that the 
restriction of the meromorphic section above to the {\sl generic fiber} $S_w$ is holomorphic, and it is the only section of the bundle $m_0K_{S/W^\prime}+ 
L$
up to a multiple, by theorem B.2.1. The term ``generic" in the sentence above means 
$$w\in W_1:= W^\prime\setminus p (R_{|S}\cup \Xi^v)$$
In the above relation, we denote by $\Xi^v$ the vertical part (with respect to $p$) of the effective divisor $\Xi_{|S}$.
By the discussion at the beginning of the paragraph A.2, the local weight of the 
metric on $m_0K_{S/W^\prime}+ L$ is given in terms of the section 
$u_{V}$ normalized in a
correct manner, as follows: for $w\in W_1$, let $\tau_w$ be the positive real number such that 
$$\tau_w^2\int_{S_w}{{\vert u_{V}\vert ^{2/m_0}}\over {|\Lambda ^rdp|^2}}
\exp(-\varphi_{\Delta^\prime})d\lambda_{S_w}= 1.\leqno (27)$$
A-priori we need to fix metrics on $S$ and $W^\prime$ in order to 
measure the norms above, but they are completely irrelevant, provided that we identify 
$\displaystyle u_{V|S_w}$ with a section of $m_0K_{S_w}+ L_{|S_w}$, therefore we skip this point, hoping that the confusion caused by it is not too big.

We denote by $\varphi^{(m_0)}_{S/W^\prime}$ the local weight of the $m_0$-Bergman kernel metric on the bundle
$m_0K_{S/W}+L$ and  then we have 
$$\exp\big(\varphi^{(m_0)}_{S/W}(x)\big)= \tau_{y}^{2m_0}|f_{V}(x)|^2\leqno(28)$$
where $y= p(x)$ and $f_{V}$ is the local expression of the section $u_{V}$ 
near the point $x$ in the fiber $S_y$. It is understood that the equality (28) holds 
for $y\in V\setminus W_1$, and that the content of the theorem 0.1 is that
the expression on the right hand side of (22) is uniformly bounded from above.

\medskip
By definition of the current $\Theta_{W^\prime}$ in (24), the exponential of the $m_0$ times its local 
potential satisfies the following relation
$$\exp\big(m_0\varphi_{\Theta_{W^\prime}}\circ p\big)= \exp\big(\varphi^{(m_0)}_{S/W^\prime}\big)
{{|f_{m_0R_{|S}}|^2}\over {|f_{m_0\Xi_{2|S}}|^2}}|f_{m_0D}|^2
\leqno (29)$$
(cf. (23) and (24) a few lines above). 
On the other hand, by relations (26) and (28) we have the equality
$$\exp\big(\varphi^{(m_0)}_{S/W^\prime}\big)
{{|f_{m_0R_{|S}}|^2}\over {|f_{m_0\Xi_{2|S}}|^2}}|f_{m_0D}|^2=
\tau_{p(x)}^{2m_0}|f_{m_0\Xi_{1|S}}|^2|f_{m_0D}|^2
\leqno (30)$$
We observe next that we have $|f_{m_0\Xi_{1|S}}|^2|f_{m_0D}|^2= |f_{m_0\Lambda}\circ p|^2$ by definition of the divisor $D$; hence we get
$$\exp\big(m_0\varphi_{\Theta_{W^\prime}}\circ p\big)=  \tau^{2m_0}|f_{m_0\Lambda}\circ p|^2\leqno (31)$$
and
it is this expression which we will use in order to evaluate the critical exponent of 
$\Theta_{W^\prime}-\Lambda$.

\smallskip
Let $\Omega\subset W^\prime$ be a coordinate open set which does not intersect the 
$p$-direct image of the divisor $R_{|S}$;
we will show next that the integral
$$\int _{w\in \Omega}{{d\lambda (w)}\over {\tau^{2+2\varepsilon_0}_w}}$$
converges, provided that the positive real number $\varepsilon_0$ is small enough.

\noindent To this end, we observe that we have 
$$\tau^2(w)\int_{S_w}{{\vert u_{m_0\Xi}\vert ^{2/m_0}}\over {|\Lambda ^rdp|^2}}\exp(-\varphi_{\Delta}-\varphi_R)d\lambda_{S_w}= 1$$
by the definition of the normalization factor $\tau$ in (27), and the 
expression of the section $u_V$ in (26).

Of course, it may happen that over some fiber 
$S_w$ the metric induced by $\Delta^\prime$ is identically $\infty$, or that the $m_0$ root of the section $u_{m_0\Xi}$ {\sl does not} belongs to the multiplier ideal of the restriction of the metric, but this kind of accidents can only happen for $w$ in an analytic set.
For such values of $w$ we simply assign the value 0 to $\tau(w)$.

We have 
$$\eqalign{
\int_{\Omega}{{d\lambda (w)}\over {\tau^{2+ 2\varepsilon}w)}}= & \int_{\Omega}d\lambda (w)\Big(\int_{S_w}{{\vert u_{m_0\Xi}\vert ^{2/m_0}}\over {|\Lambda ^rdp|^2}}\exp(-\varphi_{\Delta}-\varphi_R)
d\lambda_{S_w}\Big)^{1+ \varepsilon_0}\leq \cr
\leq & C\int_{\Omega}d\lambda (w)\int_{S_w}{{\vert u_{m_0\Xi}\vert ^{{{2+ 2\varepsilon_0}\over {m_0}}}}\over {|\Lambda ^rdp|^2}}
{{\exp\big(-(1+ \varepsilon_0)(\varphi_{\Delta^\prime} +\varphi_R)\big)}
\over {|\Lambda ^rdp|^{2\varepsilon_0}}}
d\lambda_{S_w}\cr 
\leq & C\int_{p^{-1}(\Omega)}\vert u_{m_0\Xi}\vert ^{{{2+ 2\varepsilon_0}\over {m_0}}}
{{\exp\big(-(1+ \varepsilon_0)(\varphi_{\Delta^\prime} +\varphi_R)\big)}
\over {|\Lambda ^rdp|^{2\varepsilon_0}}}d\lambda<\infty \cr
}$$
where the constant $C$ above bounds the volume of the fibers $S_w$, the first inequality is obtained by H\"older, and the last one is due to the fact that 
$\Delta_{|S}^\prime$ has trivial multiplier ideal sheaf. Therefore, the theorem 0.3 is proved.
\hfill \qed
\medskip

\noindent We will show next that by adding an arbitrarily small ample part to 
$K_X+ \Delta$ we can ``convert" the current $\Theta_{W^\prime}$ into a divisor, with the same integrability properties as in 0.3.  

The main tool for this is the fact that a closed positive current of 
$(1, 1)$--type
can be approximate in a very accurate way by 
a sequence of currents given by algebraic metrics. Let $A$ be any ample line bundle 
on $X$, and let $F_j$ be the family of divisors considered in the formula (3)
in the introduction. We consider the positive rational numbers $a^j$ such that
$$g^\star (A)-\sum a^jF_j$$
is ample on $W^\prime$. For $N\gg 0$ and divisible enough, we denote by 
$A^\prime_N$ the bundle $N\big(g^\star (A)-\sum_j a^jF_j\big)$ on $W^\prime$. The regularization theorem in [9] shows that we have
$$\varphi_{\Theta_{W^\prime}}\leq {1\over m}\log\big(\sum_{j=1}^{N_m}|h_j^{(m)}|^2\big)+ C\leqno (33)$$
where the holomorphic functions $h_j^{(m)}$ in the inequality above correspond to
global sections of the bundle 
$$m\big(g^\star (K_X+ \Delta_{|W})+ \Lambda-K_{W^\prime}
\big)+ A^\prime_N.\leqno(34)$$
We denote by $f_j$ the local equation of the hypersurface $F_j$, and then the inequality 
(33) implies
$$\varphi_{\Theta_{W^\prime}}+ {{N}\over {m}}\sum_pa^p\log |f_p|^{2}
\leq {1\over m}\log\big(\sum_{j=1}^{N_m}|l_j^{(m)}|^2\big)+ C\leqno (35)$$
where $l_j^{(m)}= h_j^{(m)}\prod f_p^{Na_p}$ are the local expression of sections of the holomorphic line bundle
$\displaystyle m\big(g^\star (K_X+ \Delta+ {N}/{m}A_{|W})+ \Lambda-K_{W^\prime}
\big)$. 

\noindent We denote by $\Theta_{W^\prime}^{(m)}$ the current associated to the metric induced by the sections
$(l^{(m)}_j)$, and we observe that we have the numerical identity
$$g^\star (K_X+ \Delta+ {N}/{m}A_{|W})+ \Lambda\equiv K_{W^\prime}+ \Theta_{W^\prime}^{(m)}.\leqno (36)$$
On the other hand, for $m\gg 0$ the relation 
$$e^{\varphi_\Lambda-\varphi_{\Theta_{W^\prime}^{(m)}}}\in L^{1+\varepsilon_1}$$
still holds on $W^\prime$ minus the image of the support of $R_S$, thanks to the 
inequality (35). Modulo a blow-up of $W^\prime$, we can assume that $\Theta_{W^\prime}^{(m)}$ is the sum of an effective divisor and of a smooth, semi-positive $(1, 1)$ form
which morally should correspond to the decomposition ``effective + nef " in the article 
[23]. \hfill\qed

\medskip

\claim B.2.3 Remark|{\rm The theorem 0.3 admits an immediate generalization 
as follows. Instead of the {\sl Weil divisor} $\Delta$ we can work with a {\sl closed positive current
with analytic singularities}. The result obtained is absolutely the same, and the modifications one has to operate in the proof presented above are minimal, so we leave them to the  interested reader. \hfill\qed
}
\endclaim

\medskip

\noindent 
\vskip 15pt

\section{References}

\bigskip

{\eightpoint

\bibitem [1]&Ambro, F.:&\ Basic properties of log canonical centers:&arXiv:math/0611205&

\bibitem [2]&Berndtsson, B.:& On the Ohsawa-Takegoshi extension theorem;& Ann.\ Inst.\ Fourier (1996)&

\bibitem [3]&Berndtsson, B.:& Integral formulas and the Ohsawa-Takegoshi extension theorem;& Science in Chi\-na, Ser A Mathematics,  2005, Vol 48&

\bibitem [4]&Berndtsson, B.:& Curvature of Vector bundles associated to holomorphic fibrations;& to appear in Ann.\ of Maths.\ (2007)&

\bibitem [5]&Berndtsson, B.:&\ An introduction to things $\dbar$;&\ to appear in the volume of the 2008 Park City Summer School, PCMI&

\bibitem [6]&Berndtsson, B., P\u aun, M.:& Bergman kernels and the pseudo-effectivity of the relative canonical bundles;& arXiv:math/0703344, Duke Math.\ Journal\ &

\bibitem [7]&Berndtsson, B., P\u aun, M.:& Qualitative extensions of twisted pluricanonical forms and closed positive currents;& on the web&

\bibitem [8]&Campana, F.:& Special Varieties and Classification Theory;&\ Annales de l'Institut Fourier 54, 3 (2004), 499--665&

\bibitem [9]&Demailly, J.-P.:&\ Regularization of closed positive currents and Intersection Theory; &J. Alg. Geom. 1 (1992) 361-409&

\bibitem [10]&Demailly, J.-P.:&  On the Ohsawa-Takegoshi-Manivel  
extension theorem;& Proceedings of the Conference in honour of the 85th birthday of Pierre Lelong, 
Paris, September 1997, Progress in Mathematics, Birkauser, 1999&

\bibitem [11]&Demailly, J.-P.:& Structure theorems for compact K\"ahler manifolds; &to appear in the volume of the 2008 Park City Summer School, PCMI&

\bibitem [12]&Ein, L., Popa, M.:&\ Extension of sections via adjoint ideals,&
arXiv:0811.4290 and
private communication, june 2007&
 
\bibitem [13]&Fujino, O.:& Higher direct images of log canonical divisors;& J. Diff.\ Geom. 66 (2004), no. 3, 453Ð479&

\bibitem [14]&Fujita, T.:&\ On Kahler fiber spaces over curves;&\ J. Math. Soc. Japan 30 (1978), no. 4, 779Ð794&

\bibitem [15]&Griffiths, P.:&\ Periods of integrals on algebraic manifolds. III. Some 
global differential-geometric properties of the period mapping;&\ Inst. Hautes 
Etudes Sci. Publ. Math. No. 38 1970 125Ð180&
 
\bibitem [16]&H\"oring, A.:& Positivity of direct image sheaves - a geometric point of view,& in preparation, available on the author's home page&

\bibitem [17]&Kawamata, Y.:&\ Pluricanonical systems on minimal algebraic varieties;& Invent. Math.  79  (1985),  no. 3&

\bibitem [18]&Kawamata, Y.:&\ On Fujita's freeness conjecture for
3-folds and 4-folds&  Math.\ Ann.\ {\bf 308}, 1997&

\bibitem [19]&Kawamata, Y.:&\ Subadjunction of log canonical divisors, 2;& Amer.\ J.\ Math.\  {\bf 120} (1998) 893--899&

\bibitem [20]&Kawamata, Y.:&\ On algebraic fiber spaces;& math.AG/0107160. in 
Contemporary Trends in Algebraic Geometry and Algebraic Topology, World 
ScientiÞc, 2002, 135--154& 

\bibitem [21]&Kawamata, Y.:&\ Semipositivity theorem for reducible algebraic fiber spaces;& rXiv:0911.1670&

\bibitem [22]&Koll\'ar, J.:&\ Higher direct images of dualizing sheaves, II;&\ Ann. of 
Math. (2) 124 (1986), no. 1, 171--202&

\bibitem [23]&Koll\'ar, J.:&\ Kodaira's canonical bundle formula and adjunction;&\ in {\sl Flips for 3-folds and 4-folds}, Oxford University Press, Oxford, 2007&

 \bibitem [24]&Manivel, L.:& Un th\'eor\`eme de prolongement L2 de sections holomorphes d'un 
 fibr\'e hermitien;& Math.\ Zeitschrift {\bf 212} (1993), 107-122&

 \bibitem [25]&McNeal J., Varolin D.:&\ Analytic inversion of adjunction: $L\sp
2$ extension theorems with gain;& Ann. Inst. Fourier (Grenoble) 57
(2007), no. 3, 703--718& 

\bibitem [26]&Mourougane, C., Takayama, S.:&\ Hodge metrics and the curvature of higher direct images;& arXiv:0707.3551&

\bibitem [27]&Mourougane, C., Takayama, S.:&\ Extension of twisted Hodge metrics for K\"ahler morphisms;& arXiv: 0809.3221&
 
\bibitem [28]&Nakayama, N.:&\ Hodge filtrations and the higher direct images of 
canonical sheaves;&\ Invent. Math. 85 (1986), no. 1, 217--221&
  
 \bibitem [29]&Ohsawa, T., Takegoshi, K.\ :& On the extension of $L^2$
holomorphic functions;& Math.\ Z.,
{\bf 195} (1987), 197--204&
 
\bibitem [30]&Tsuji, H.:& Extension of log pluricanonical forms from subvarieties;& math.CV/0511342 &

\bibitem [31]&Schmid, W.:&\ Variation of Hodge structure: the singularities of the 
period mapping;&\ Invent. Math. 22 (1973), 211Ð319&
   
\bibitem [32]&Siu, Y.-T.:& Invariance of Plurigenera;& Inv.\ Math.,
{\bf 134} (1998), 661-673&

\bibitem [33]&Siu, Y.-T.:& Extension of twisted pluricanonical sections with plurisubharmonic weight and invariance of semipositively twisted plurigenera for manifolds not necessarily of general type;& Complex geometry (G\"ottingen, {\bf 2000}),  223--277, Springer, Berlin, 2002&

\bibitem [33]&Siu, Y.-T.:& Multiplier ideal sheaves in complex and algebraic geometry;& Sci.\ China Ser.  {\bf A 48}, 2005&

\bibitem [34]&Viehweg, E.:&\ Weak positivity and the additivity of the Kodaira dimension for certain fibre spaces. 
Proc. Algebraic Varieties and Analytic Varieties,&\ Adv. Studies in 
Math. 1, Kinokunya--North-Holland Publ. 1983, 329--353& 
     
\bibitem [35]&Viehweg, E.:&\ Quasi-Projective
Moduli for
Polarized Manifolds;& Springer-Verlag, Berlin, Heidelberg, New York, 1995
as: Ergebnisse der Mathematik und ihrer Grenzgebiete, 3. Folge, Band 30&

}

\bigskip
\noindent
(version of February 20, 2010, printed on \today)
\bigskip\bigskip
{\parindent=0cm
Bo Berndtsson,  
bob@math.chalmers.se\\
Mihai P\u aun,
paun@iecn.u-nancy.fr
}

\end